\documentclass{amsart}
\usepackage{amssymb}

\def\RR{\mathbb R}

\def\gam{\gamma }
\def\lam{\lambda }
\def\eps{\varepsilon }

\newcommand{\remove}[1]{{}}
\newcommand{\set}[1]{\{ #1 \}}

\newtheorem{theorem}{Theorem}[section]

\newtheorem{lemma}[theorem]{Lemma}
\newtheorem{corollary}[theorem]{Corollary}

\theoremstyle{definition}
\newtheorem*{definition}{Definition}

\theoremstyle{remark}
\newtheorem{remark}{Remark}[section]

\newtheorem{example}{Example}[section]
\newtheorem{examples}[example]{Examples}

\numberwithin{equation}{section}

\title[Ternary alphabets]{Generalized golden ratios of ternary alphabets}

\author{Vilmos Komornik}
\address{Delft Institute of Applied Mathematics\\
Technical University of Delft\\
Mekelweg 4\\
2628 CD Delft\\
The Netherlands\\
and 
D\'epartement de math\'ematique\\
         Universit\'e de Strasbourg\\
         7 rue Ren\'e Descartes\\
         67084 Strasbourg Cedex, France}
\email{komornik@math.u-strasbg.fr}
\author{Anna Chiara Lai}
\address{Dipartimento Me.Mo.Mat.\\
          Universit\`a degli Studi di Roma ``La Sapienza''\\
          Via A. Scarpa, 16\\
          00161 Roma, Italy}
\email{lai@dmmm.uniroma1.it}
\author{Marco Pedicini}
\address{Istituto per le Applicazioni del Calcolo ``Mauro Picone''\\
          Consiglio Nazionale delle Ricerche\\
          Viale del Policlinico, 137\\
          00161 Roma, Italy}
\email{marco@iac.cnr.it}
\subjclass[2000]{Primary: 11A63, Secondary: 11B83}
\keywords{Golden ratio, ternary alphabet, unique expansion, noninteger base, beta-expansion, greedy expansion, lazy expansion,
univoque sequence}

\begin{document}

\begin{abstract}
Expansions in noninteger bases often appear in number theory and probability theory, and they are closely connected to ergodic theory, measure theory and topology. For two-letter alphabets the golden ratio plays a special role: in smaller bases only trivial expansions are unique, whereas in greater bases there exist nontrivial unique expansions. In this paper we determine the corresponding critical bases for all three-letter alphabets and we establish the fractal nature of these bases in function of the alphabets. 
\end{abstract}

\maketitle

\section{Introduction}\label{s1}

Since the appearence of R\'enyi's $\beta$-expansions \cite{R} many works were devoted to expansions in noninteger bases. Much research was stimulated by the discovery of Erd\H os, Horv\'ath and Jo\'o \cite{EHJ} who proved the existence of many real numbers $1<q<2$ for which only one sequence $(c_i)$ of zeroes and ones satisfies the equality
\begin{equation*}
\sum_{i=1}^{\infty}\frac{c_i}{q^i}=1.
\end{equation*}
The set of such ``univoque'' bases has a fractal nature, see, e.g., \cite{EJK1}, \cite{KL1}, \cite{KL3}, where arbitrary bases $q>1$ are also considered. 

Contrary to the integer case, in a given noninteger base $q>1$ a real number $x$ may have sometimes many different expansions of the form
\begin{equation}\label{11}
\pi_q(c):=\sum_{i=1}^{\infty}\frac{c_i}{q^i}=x
\end{equation}
with integer ``digits'' satisfying $0\le c_i<q$ for every $i$. On the other hand, the set of numbers $x$ having a unique expansion has many unexpected topological and combinatorial properties, depending on the value of $q$; see, e.g., Dar\'oczy and K\'atai \cite{DK1}, Glendinning and Sidorov \cite{GS}, and \cite{DeVK}.

Given a finite alphabet $A=\set{a_1,\ldots, a_J}$ of real numbers $a_1<\cdots <a_J$ and a real number $q>1$, by an expansion of a real number $x$ we mean a sequence $(c_i)$ of numbers $c_i\in A$ satisfying \eqref{11}. The expansions of
\begin{equation*}
x_1:=\sum_{i=1}^{\infty}\frac{a_1}{q^i}\quad\text{and}\quad x_2:=\sum_{i=1}^{\infty}\frac{a_J}{q^i}
\end{equation*}
are always unique; they are called the trivial unique expansions.

For two-letter alphabets $A=\set{a_1, a_2}$ the golden ratio $p:=(1+\sqrt{5})/2$ plays a special role: there exist nontrivial unique expansions in base $q$ if and only if $q>p$. 

The purpose of this paper is to determine analogous critical bases for each ternary alphabet $A=\set{a_1, a_2, a_3}$. Our main tool is a lexicographic characterization of unique expansions, given in \cite{Ped05}, which generalized for arbitrary finite alphabets a theorem of Parry \cite{P} and its various extensions  \cite{DK1}, \cite{EHJ}, \cite{EJK1}, \cite{KL2}.

By a normalization it suffices to consider the alphabets $A_m:=\set{0,1,m}$ with $m\ge 2$. Our main result is the following:

\begin{theorem}\label{t11}\mbox{}
There exists a continuous function $p:[2,\infty)\to\RR$, $m\mapsto p_m$ satisfying
\begin{equation*}
2\le p_m\le P_m:=1+\sqrt{\frac{m}{m-1}}
\end{equation*}
for all $m$ such that the following properties hold true:
\smallskip

(a) for each $m\ge 2$, there exist nontrivial unique expansions if and only if $q>p_m$;

\smallskip
(b) we have $p_m=2$ if and only if $m=2^k$ for some positive integer $k$;
\smallskip

(c) the set $C:=\set{m\ge 2\ :\ p_m= P_m}$ is a Cantor set, i.e., an uncountable closed set having neither interior nor isolated points; its smallest element is $1+x\approx 2.3247$ where $x$ is the first Pisot number, i.e., the positive root of the equation $x^3=x+1$;
\smallskip

(d) each connected component $(m_d,M_d)$ of $[2,\infty)\setminus C$ has a point $\mu_d$ such that $p$ is strictly decreasing in $[m_d,\mu_d]$ and strictly increasing in $[\mu_d,M_d]$.
\end{theorem}

Moreover, we will determine explicitly the function $p$ and the  numbers $m_d$, $M_d$, $\mu_d$.

In Section \ref{s2}, we review some basic facts about expansions and we also give some new results. In Sections \ref{s3}--\ref{s4} we introduce the class of \emph{admissible sequences} and we clarify their structure and their basic properties. These results allow us  to determine in Section \ref{s5} the critical bases for all ternary alphabets. 

\section{Some results on arbitrary alphabets}\label{s2}

Throughout this section we consider a fixed finite alphabet $A=\set{a_1,\ldots, a_J}$ of real numbers $a_1<\cdots <a_J$. Given a real number $q>1$, by an expansion of a real number $x$ we mean a sequence $(c_i)$ of numbers $c_i\in A$ satisfying the equality
\begin{equation*}
\pi_q(c):=\sum_{i=1}^{\infty}\frac{c_i}{q^i}=x.
\end{equation*}
In order to have an expansion, $x$ must belong to the interval $[\frac{a_1}{q-1},\frac{a_J}{q-1}]$. Conversely, we recall from \cite{Ped05} the following results:

\begin{theorem}\label{t21}
Every $x\in [\frac{a_1}{q-1},\frac{a_J}{q-1}]$ has at least one expansion in base $q$ if and only if
\begin{equation}\label{21}
1<q\le Q_A:=1+\frac{a_J-a_1}{\max_{j>1} \set{a_j-a_{j-1}}} (\le J).
\end{equation}
\end{theorem}

A sequence $(c_i)\in A^{\infty}$ is called \emph{univoque} in base $q$ if
\begin{equation*}
x:=\sum_{i=1}^{\infty}\frac{c_i}{q^i}
\end{equation*}
has no other expansion in this base. The constant sequences $(a_1)^{\infty}$ and $(a_J)^{\infty}$ are univoque in every base $q$: they are called the \emph{trivial unique expansions}. We also recall from \cite{Ped05} the following characterization of unique expansions:

\begin{theorem}\label{t22}
Assume \eqref{21}. An expansion $(c_i)$ is unique in base $q$ if and only if the following conditions are satisfied:
\begin{align*}
&\sum_{i=1}^\infty \frac{c_{n+i}-a_1}{q^i}  < a_{j+1}-a_j \quad\text{whenever}\quad  c_n=a_j<a_J;\\
&\sum_{i=1}^\infty \frac{a_J-c_{n+i}}{q^i}  < a_j-a_{j-1} \quad\text{whenever}\quad  c_n=a_j>a_1.
\end{align*}
\end{theorem}

\begin{proof}[Proof of the sufficiency]
We have to show that if $(d_i)$ is another sequence in $A$ then it represents a different sum. Let $n\ge 1$ be the first index such that $c_n\ne d_n$. If $c_n<d_n$, then writing $c_n=a_j$ we have $a_j<a_J$, so that
\begin{equation*}
\sum_{i=1}^\infty \frac{d_i}{q^i}-\sum_{i=1}^\infty \frac{c_i}{q^i}
\ge \frac{a_{j+1}-a_j}{q^n}-\sum_{i=n+1}^\infty\frac{a_1-c_i}{q^i}>0
\end{equation*}
by our assumption. If $c_n>d_n$, then writing $c_n=a_j$ we have $a_j>a_1$, so that
\begin{equation*}
\sum_{i=1}^\infty \frac{c_i}{q^i}-\sum_{i=1}^\infty \frac{d_i}{q^i}
\ge \frac{a_j-a_{j-1}}{q^n}+\sum_{i=n+1}^\infty\frac{c_i-a_J}{q^i}>0
\end{equation*}
by our second assumption.
\end{proof}

\begin{proof}[Proof of the necessity]
If the first condition is not satisfied for some $c_n=a_j<a_J$, then by Theorem \ref{t21} there exists another expansion beginning with $c_1\cdots c_{n-1}a_{j+1}$. If the second condition is not satisfied for some $c_n=a_j>a_1$, then by Theorem \ref{t21} there exists another expansion beginning with $c_1\cdots c_{n-1}a_{j-1}$.
\end{proof}

Let us mention some consequences of this characterization.

\begin{corollary}\label{c23}
For every given set  $C\subset A^{\infty}$ there exists a number
\begin{equation*}
1\le q_C\le Q_A
\end{equation*}
such that
\begin{align*}
q>q_C&\Longrightarrow \text{ every sequence $c\in C$ is univoque in base }q;\\
1<q<q_C&\Longrightarrow\ \text{ not every sequence $c\in C$ is univoque in base }q.
\end{align*}
\end{corollary}

\begin{proof}
If $C=\varnothing$, then we may choose $q_C=1$. If $C$ is nonempty, then for each sequence $c\in C$, each condition of Theorem \ref{t22} is equivalent to an inequality of the form $q>q_{\alpha}$. Since we consider only bases $q$ satisfying \eqref{21}, we may assume that
$q_{\alpha}\le Q_A$ for every $\alpha$. Then
\begin{equation*}
q_C:=\max\set{1,\sup q_{\alpha}}
\end{equation*}
 has the required properties.
\end{proof}

\begin{definition}\label{d21}
The number $q_C$ is called the \emph{critical base} of $C$. If $C=\set{c}$ is a one-point set, then $q_c:=q_C$ is also called the critical base of the sequence $c$.
\end{definition}

\begin{remark}\label{r21}
If $C$ is a nonempty finite set of eventually periodic sequences, then the supremum $\sup q_{\alpha}$ in the above proof is actually a maximum. In this case not all sequences $c\in C$ are univoque in base $q=q_C$.
\end{remark}

\begin{example}\label{e21}
Consider the ternary alphabet $A=\set{0,1,3}$ and the periodic sequence $(c_i)=(31)^\infty$. By the periodicity of $(c_i)$ we have for each $n$ either $c_n=3$ and $(c_{n+i})=(13)^\infty$ or $c_n=1$ and $(c_{n+i})=(31)^\infty$. According to the preceding remark Theorem \ref{t22} contains only three conditions on $q$.  For $c_n=3$ we have the condition
\begin{equation*}
\sum_{i=1}^\infty \frac{3-c_{n+i}}{q^i}<2\Longleftrightarrow \frac{2q}{q^2-1}<2,
\end{equation*}
while for $c_n=3$ we have the following two conditions:
\begin{equation*}
\sum_{i=1}^\infty \frac{3-c_{n+i}}{q^i}<1\Longleftrightarrow \frac{2}{q^2-1}<1
\end{equation*}
and
\begin{equation*}
\sum_{i=1}^\infty \frac{c_{n+i}}{q^i}<2\Longleftrightarrow \frac{3}{q-1}-\frac{2}{q^2-1}<2.
\end{equation*}
They are equivalent approximatively to the inequalities $q>1.61803$, $q>1.73205$ and $q>2.18614$ respectively, so that $q_c\approx 2.18614$.
\end{example}

It is well-known that for the alphabet $A=\set{0,1}$ there exist nontrivial univoque sequences in base $q$ if and only if $q>\frac{1+\sqrt{5}}{2}$.
There exists a ``generalized golden ratio'' for every alphabet:

\begin{corollary}\label{c24}
There exists a number $1< G_A\le Q_A$ such that
\begin{align*}
q>G_A&\Longrightarrow \text{ there exist nontrivial univoque sequences};\\
1<q<G_A&\Longrightarrow\ \text{ there are no nontrivial univoque sequences}.
\end{align*}
\end{corollary}

\begin{proof}
If a sequence is univoque in some base, then it is also univoque in every larger base. If there exists a base satisfying \eqref{21} in which there exist nontrivial univoque sequences, then it follows that the infimum of such bases satisfies the requirements for $G_A$, except perhaps the strict inequality $G_A>1$. Otherwise we may choose $G_A:=Q_A$.

It remains to prove that if $q>1$ is sufficiently close to one, then the only univoque sequences are $a_1^{\infty}$ and $a_J^{\infty}$. We show that it suffices to choose $q>1$ so small that the following three conditions are satisfied:
\begin{align}
&\frac{a_J-a_1}{q-1}\ge a_{j+1}-a_{j-1},\quad j=2,\ldots, J-1,\label{22}\\
&\frac{a_j-a_1}{q}+\frac{1}{q}\cdot\frac{a_J-a_1}{q-1}\ge (a_2-a_1)+\frac{a_j-a_{j-1}}{q},\quad j=2,\ldots, J,\label{23}\\
&\frac{a_J-a_j}{q}+\frac{1}{q}\cdot\frac{a_J-a_1}{q-1}\ge (a_J-a_1)+\frac{a_{j+1}-a_j}{q},\quad j=1,\ldots, J-1.\label{24}
\end{align}
The proof consists of three steps. Let $(c_i)$ be a univoque sequence in base $q$.

If $c_n=a_j$ for some $n$ and $1<j<J$, then the conditions of Theorem \ref{t22} imply that
\begin{equation*}
\sum_{i=1}^\infty \frac{c_{n+i}-a_1}{q^i}  < a_{j+1}-a_j
\quad\text{and}\quad
\sum_{i=1}^\infty \frac{a_J-c_{n+i}}{q^i}  < a_j-a_{j-1}.
\end{equation*}
Taking their sum we conclude that
\begin{equation*}
\frac{a_J-a_1}{q-1}<a_{j+1}-a_{j-1},
\end{equation*}
which contradicts \eqref{22}. This proves that $c_n\in\set{a_1,a_J}$ for every $n$.

If $c_n=a_1$ and $c_{n+1}=a_j>a_1$ for some $n$, then applying Theorem \ref{t22} we obtain that
\begin{equation*}
\sum_{i=1}^\infty \frac{c_{n+i}-a_1}{q^i}  < a_2-a_1
\quad\text{and}\quad
\sum_{i=1}^\infty \frac{a_J-c_{n+i+1}}{q^i}  < a_j-a_{j-1}.
\end{equation*}
Dividing the second inequality by $q$ and adding the result to the first one we obtain that
\begin{equation*}
\frac{a_j-a_1}{q}+\frac{1}{q}\cdot\frac{a_J-a_1}{q-1}<(a_2-a_1)+\frac{a_j-a_{j-1}}{q},
\end{equation*}
which contradicts \eqref{23}. This proves that $c_n=a_1$ implies $c_{n+1}=a_1$ for every $n$.

Finally, if $c_n=a_J$ and $c_{n+1}=a_j<a_J$ for some $n$, then applying Theorem \ref{t22} we obtain that
\begin{equation*}
\sum_{i=1}^\infty \frac{a_J-c_{n+i}}{q^i}  < a_J-a_{J-1}
\quad\text{and}\quad
\sum_{i=1}^\infty \frac{c_{n+i+1}-a_1}{q^i}  < a_{j+1}-a_j.
\end{equation*}
Dividing the second inequality by $q$ and adding the result to the first one we now obtain that
\begin{equation*}
\frac{a_J-a_j}{q}+\frac{1}{q}\cdot\frac{a_J-a_1}{q-1}<(a_J-a_{J-1})+\frac{a_{j+1}-a_j}{q},
\end{equation*}
which contradicts \eqref{24}. This proves that $c_n=a_J$ implies $c_{n+1}=a_J$ for every $n$.
\end{proof}

\begin{definition}\label{d22}
The number $G_A$ is called the \emph{critical base} of the alphabet $A$.
\end{definition}

The following invariance properties of critical bases readily follow from the definitions; they will simplify our proofs.

\begin{lemma}\label{l25}
The critical base does not change if we replace the alphabet $A$
\begin{itemize}
\item by $b+A=\set{b+a_j\ |\ j=1,\ldots,m}$ for some real number $b$;
\item by $dA=\set{da_j\ |\ j=1,\ldots,m}$ for some nonzero real number $d$;
\item by the \emph{conjugate alphabet} $A':=\set{a_m+a_1-a_j\ |\ j=1,\ldots,m}$.
\end{itemize}
\end{lemma}

\begin{proof}
First we note that $Q_A=Q_{b+A}=Q_{dA}=Q_{A'}$. Fix a base $1<q\le Q_A$ and a sequence $(c_i)$ of real numbers.
It follows from the definitions that the following properties are equivalent:

\begin{itemize}
\item $(c_i)$ is an expansion of $x$ for the alphabet $A$;

\item $(b+c_i)$ is an expansion of $x+\frac{b}{q-1}$ for the alphabet $b+A$;

\item $(dc_i)$ is an expansion of $dx$ for the alphabet $dA$;

\item $(a_m+a_1-c_i)$ is an expansion of $\frac{a_m+a_1}{q-1}-x$ for the alphabet $A'$.
\end{itemize}
Hence if one of these expansions is unique, then the others are unique as well.
\end{proof}

\section{Admissible sequences}
\label{s3}

This section contains some preliminary technical results.

\begin{definition}
A sequence $d=(d_i)=d_1d_2\cdots$ of zeroes and ones is \emph{admissible} if
\begin{equation}\label{31}
0d_2d_3\cdots \le (d_{n+i})\le d_1d_2d_3\cdots
\end{equation}
for all $n=0,1,\ldots .$
\end{definition}

\begin{examples}\mbox{}\label{e31}

\begin{itemize}
\item The trivial sequences  $0^{\infty}$ and $1^{\infty}$ are admissible.

\item More generally, the sequences $(1^N0)^{\infty}$ ($N=1,2,\ldots$) and  $(10^N)^{\infty}$ ($N=0,1,\ldots$) are admissible.

\item The sequence $(11010)^{\infty}$ is also admissible.

\item The (not purely periodic) sequence $10^{\infty}$ is admissible.
\end{itemize}
\end{examples}

In order to clarify the structure of admissible sequences we give an equivalent recursive definition. Given a sequence $h=(h_i)$ of positive integers, starting with
\begin{equation*}
S_h(0,1):=1\quad\text{and}\quad S_h(0,0):=0
\end{equation*}
we define the blocks $S_h(j,1)$ and $S_h(j,0)$ for $j=1,2,\ldots$ by the recursive formulae
\begin{align*}
&S_h(j,1):=S_h(j-1,1)^{h_j}S_h(j-1,0)
\intertext{and}
&S_h(j,0):=S_h(j-1,1)^{h_j-1}S_h(j-1,0).
\end{align*}
Observe that $S_h(j,1)$ and $S_h(j,0)$ depend only on $h_1$,\ldots, $h_j$, so that they can also be defined for every finite sequence $h=(h_j)$ of length $\ge j$.
We also note that $S_h(j,0)=S_h(j-1,0)$ whenever $h_j=1$ and that the length $\ell_j$ of $S_h(j,1)$ tends to infinity as $j\to\infty$.

If the sequence $h=(h_j)$ is given, we often omit the subscript $h$ and we write simply $S(j,1)$ and $S(j,0)$.

Let us mention some properties of these blocks that we use in the sequel. Given two finite blocks $A$ and $B$ we write for brevity

\begin{itemize}
\item $A\to B$ or $B=\cdots A$ if $B$ ends with $A$;

\item $A<B$ or $A\cdots <B\cdots $ if $Aa_1a_2\cdots <Bb_1b_2\cdots $ lexicographically for any sequences $(a_i)$ and $(b_i)$ of zeroes and ones;

\item $A\le B$ or $A\cdots \le B\cdots $ if $A<B$ or $A=B$.
\end{itemize}

\begin{lemma}\label{l31}
For any given sequence $h=(h_j)$ the blocks $S(j,1)$ and $S(j,0)$ have the following properties:\mbox{}

(a) We have
\begin{align}
&S(j,1)=1S(1,0)\cdots S(j,0)\quad\text{for all}\quad j\ge 0;\label{32}\\
&S(0,0)\cdots S(j-1,0)\to S(j,1)\quad\text{for all}\quad j\ge 1;\label{33}\\
&S(0,0)\cdots S(j-1,0)\to S(j,0)\quad\text{whenever}\quad h_j\ge 2;\label{34}\\
&S(j,0)<S(j,1)\quad\text{for all}\quad j\ge 0.\label{35}
\end{align}
\smallskip

(b) If $A_j\to S(j,1)$ for some nonempty block $A_j$, then $A_j\le S(j,1)$.
\smallskip

(c) If $B_j\to S(j,0)$ for some nonempty block $B_j$, then $B_j\le S(j,0)$.
\end{lemma}

\begin{proof}\mbox{}

(a) \emph{Proof of \eqref{32}.} For $j=0$ we have $S(j,1)=1$ by definition. If $j\ge 1$ and the identity is true for $j-1$, then the identity for $j$ follows by using the equality $S(j,1)=S(j-1,1)S(j,0)$ coming from the definition of $S(j,1)$ and $S(j-1,1)$.
\smallskip

\emph{Proof of \eqref{33} and \eqref{34}.} For $j=1$ we have $S(0,0)=0$ and $S(1,0)=1^{h_1-1}0$, so that $S(0,0)\to S(1,0)\to S(1,1)$. (The condition $h_1\ge 2$ is not needed here.) Proceeding by induction, if \eqref{33} holds for some $j\ge 1$, then both hold for $j+1$ because
\begin{equation*}
S(0,0)\cdots S(j-1,0)S(j,0)\to S(j,1)S(j,0)\to S(j+1,1),
\end{equation*}
and in case $h_{j+1}\ge 2$ we have also $S(j,1)S(j,0)\to S(j+1,0)$.
\smallskip

\emph{Proof of \eqref{35}.} The case $j=0$ is obvious because the left side begins with $0$ and the right side begins with $1$. If $j\ge 1$ and \eqref{35} holds for $j-1$, then we deduce from the inequality $S(j-1,0)\cdots <S(j-1,1)\cdots$ that
\begin{equation*}
S(j,0)\cdots =S(j-1,1)^{h_j-1}S(j-1,0)\cdots<S(j-1,1)^{h_j}\cdots .
\end{equation*}
Since $S(j,1)$ begins with $S(j-1,1)^{h_j}$, this implies \eqref{35} for $j$.
\smallskip

(b) We may assume that $A_j\ne S(j,1)$; this excludes the case $j=0$ when we have necessarily $A_0=S(0,1)=1$. For $j=1$ we have $S(j,1)=1^{h_1}0$ and $A_j=1^t0$ with some integer $0\le t<h_1$, and we conclude by observing that $1^t0\cdots <1^{h_1}\cdots .$

Now let $j\ge 2$ and assume that the result holds for $j-1$. Using the equality $S(j,1)=S(j-1,1)^{h_j}S(j-1,0)$ we distinguish three cases.

If $A_j\to S(j-1,0)$, then we have the implications
\begin{align*}
A_j\to S(j-1,0)
&\Longrightarrow A_j\to S(j-1,1)\quad\text{and}\quad A_j\ne S(j-1,1)\\
&\Longrightarrow A_j\cdots <S(j-1,1)\cdots \\
&\Longrightarrow A_j\cdots <S(j,1)\cdots .
\end{align*}

If $A_j=A_{j-1}S(j-1,1)^tS(j-1,0)$ for some $0\le t<h_j$, $A_{j-1}\to S(j-1,1)$ and $A_{j-1}\ne S(j-1,1)$, then
\begin{align*}
A_{j-1}\cdots < S(j-1,1)\cdots
&\Longrightarrow A_j\cdots <S(j-1,1)\cdots \\
&\Longrightarrow A_j\cdots <S(j,1)\cdots .
\end{align*}

Finally, if $A_j=S(j-1,1)^tS(j-1,0)$ for some $0\le t<h_j$, then using \eqref{35} we have
\begin{equation*}
A_j\cdots <S(j-1,1)^{t+1}\cdots
\end{equation*}
and therefore
\begin{equation*}
A_j\cdots <S(j,1)\cdots .
\end{equation*}
\smallskip

(c) Proceeding by induction, the case $j=0$ is obvious because then we have necessarily $B_0=S(0,0)=0$. Let $j\ge 1$ and assume that the property holds for $j-1$ instead of $j$. If $h_j>1$, then the case of $j$ follows by applying part (b) with $h_j$ replaced by $h_j-1$. If $h_j=1$, then we have $S(j,0)=S(j-1,0)$ and applying (b) we conclude that
\begin{equation*}
B_j\to S(j,0)\Longrightarrow B_j\to S(j-1,0)\Longrightarrow B_j\le S(j-1,0)=S(j,0).\qedhere
\end{equation*}
\end{proof}

The following lemma is a partial converse of \eqref{33}.

\begin{lemma}\label{l32}
If $A$ is a block of length $\ell_{N-1}$ in some sequence $S(N,a_1)S(N,a_2)\cdots$ with $N\ge 1$ and $(a_i)\subset\set{0,1}$, then $A\ge S(0,0)\cdots S(N-1,0)$. Furthermore, we have $A= S(0,0)\cdots S(N-1,0)$ if and only if $A\to S(N,a_i)$ with some $a_i=1$.
\end{lemma}

\begin{proof}
The case $N=1$ is obvious because then $S(0,0)=0$ implies that $A=0$, and $S(1,1)=1^{h_1}0$ ends with $0$.

Now let $N\ge 2$ and assume by induction that the result holds for $N-1$. Writing $A=BC$ with a block $B$ of the same length as $S(0,0)\cdots S(N-2,0)$ and applying the induction hypothesis to $B$ in the sequence
\begin{equation*}
S(N,a_1)S(N,a_2)\cdots= \left( S(N-1,1)^{h_N-1+a_i}S(N-1,0)\right)
\end{equation*}
we obtain that $B\to S(N-1,1)$ for one of the blocks on the right side and thus $B=S(0,0)\cdots S(N-2,0)$. Then it follows from our assumption that $C$ has the same length as $S(N-1,0)$ and $C\le S(N-1,0)$. Since $S(N-1,0)<S(N-1,1)$, the block containing $B$ must be followed by a block $S(N-1,0)$. We conclude that $C=S(N-1,0)$ and therefore $A=BC=S(0,0)\cdots S(N-1,0)$ and
\begin{equation*}
A\to S(N-1,1)^{h_N-1+a_i}S(N-1,0)=S(N,a_i)
\end{equation*}
for some $a_i=1$.
\end{proof}

\begin{lemma}\label{l33}
A sequence $d=(d_i)$ is admissible if and only if one of the following three conditions is satisfied:

\begin{itemize}
\item $d=0^{\infty}$;

\item there exists an infinite sequence $h=(h_i)$ of positive integers such that $d$ begins with $S_h(N,1)$ for every $N=0,1,\ldots $;

\item $d=S_h(N,1)^{\infty}$ with some nonnegative integer $N$ and a finite sequence $h=(h_1,\ldots,h_N)$ of positive integers.
\end{itemize}

\end{lemma}

\begin{proof}
It follows from the definition that $d_1=1$ for all admissible sequences other than $0^{\infty}$. In the sequel we consider only admissible sequences beginning with $d_1=1$. We omit the subscript $h$ for brevity.

Let $d=(d_i)$ be an admissible sequence. Setting $d^0_i:=d_i$ we have
\begin{equation*}
d=S(0,d^0_1)S(0,d^0_2)\cdots
\end{equation*}
with the admissible sequence  $(d^0_i)$.

Proceeding by recurrence, assume that
\begin{equation*}
d=S(j,d^j_1)S(j,d^j_2)\cdots
\end{equation*}
for some integer $j\ge 0$ with an admissible sequence $(d^j_i)$ and positive integers $h_1$,\ldots, $h_j$. (We need no such positive integers for $j=0$.)

If $(d^j_i)=1^{\infty}$, then $d=S(j,1)^{\infty}$. Otherwise there exists a positive integer $h_{j+1}$ such that $d$ begins with $S(j,1)^{h_{j+1}}S(j,0)$. Since the sequence $(d^j_i)$ is admissible, we have
\begin{equation*}
S(j,0)S(j,1)^{h_{j+1}-1}S(j,0)\cdots \le S(j,d^j_{n+1})S(j,d^j_{n+2})\cdots\le S(j,1)^{h_{j+1}}S(j,0)\cdots
\end{equation*}
for all $n=0,1, \ldots .$ This implies that each block $S(j,0)$ is followed by at least $h_{j+1}-1$ and at most $h_{j+1}$ consecutive blocks $S(j,1)$, so that
\begin{equation*}
d=S(j+1,d^{j+1}_1)S(j+1,d^{j+1}_2)\cdots
\end{equation*}
for a suitable sequence $(d^{j+1}_i)$ of zeroes and ones. The admissibility of $(d^j_i)$ can then be rewritten in the form
\begin{multline}
S(j,0)S(j+1,0)S(j+1,d^{j+1}_2)S(j+1,d^{j+1}_3)\cdots\\
\le S(j,d^j_{n+1})S(j,d^j_{n+2})\cdots\\
\le S(j+1,1)S(j+1,d^{j+1}_2)S(j+1,d^{j+1}_3)\cdots\label{36}
\end{multline}
for $n=0,1, \ldots .$

We claim that the sequence $(d^{j+1}_i)$ is also admissible. We have $d^{j+1}_1=1$ by the definition of $h_{j+1}$. It remains to show that
\begin{multline*}
S(j+1,0)S(j+1,d^{j+1}_2)S(j+1,d^{j+1}_3)\cdots\\
\le S(j+1,d^{j+1}_{k+1})S(j+1,d^{j+1}_{k+2})S(j+1,d^{j+1}_{k+3})\cdots\\
\le S(j+1,1)S(j+1,d^{j+1}_2)S(j+1,d^{j+1}_3)\cdots
\end{multline*}
for $k=0,1, \ldots .$

The second inequality is a special case of the second inequality of \eqref{36}. The first inequality is obvious for $k=0$. For $k\ge 1$ it is equivalent to
\begin{multline*}
S(j,0)S(j+1,0)S(j+1,d^{j+1}_2)S(j+1,d^{j+1}_3)\cdots\\
\le S(j,0)S(j+1,d^{j+1}_{k+1})S(j+1,d^{j+1}_{k+2})S(j+1,d^{j+1}_{k+3})\cdots
\end{multline*}
and this is a special case of the first inequality of \eqref{36} because $S(j+1,d^{j+1}_k)$ ends with $S(j,0)$.

It follows from the above construction that $(d_i)$ has one of the two forms specified in the statement of the lemma.
\medskip

Turning to the proof of the converse statement, first we observe that if $d$ begins with $S(N,1)$ for some sequence $h=(h_i)$ and for some integer $N\ge 1$, then
\begin{equation}\label{37}
d_n\cdots d_{\ell_N} < d_1\cdots d_{\ell_N-n+1}\quad\text{for}\quad n=2,\ldots, \ell_N;
\end{equation}
this is just a reformulation of part (b) of Lemma \ref{l31}.

If $d_1d_2\cdots$ begins with $S(N,1)$ for all $N$, then the second inequality of \eqref{31} follows for all $n\ge 1$ by using the relation $\ell_N\to\infty$. Moreover, the inequality is strict. For $n=0$ we have clearly equality.

If $d=S(N,1)^{\infty}$ for some $N\ge 0$, then $d$ is $\ell_N$-periodical so that the second inequality of \eqref{31} follows from \eqref{37} for
all $n$, except if $n$ is a multiple of $\ell_N$; we get strict inequalities in theses cases. If $n$ is a multiple of $\ell_N$, then we have obviously equality again.

It remains to prove the first inequality of \eqref{31}. If $d=S(N,1)^{\infty}$ for some $N\ge 0$, then
we deduce from Lemma \ref{l32} that either
\begin{equation*}
(d_{n+i})>S(0,0)\cdots S(N-1,0)
\end{equation*}
or
\begin{equation*}
(d_{n+i})=S(0,0)\cdots S(N-1,0)S(N,1)^{\infty}.
\end{equation*}
Since
\begin{equation*}
0d_2d_3\cdots =S(0,0)\cdots S(N-1,0)S(N,0)S(N,1)^{\infty},
\end{equation*}
we conclude in both cases the  strict inequalities
\begin{equation*}
(d_{n+i})>0d_2d_3\cdots .
\end{equation*}

If  $d_1d_2\cdots$ begins with $S(N,1)$ for all $N$, then
\begin{equation*}
0d_2d_3\cdots =S(0,0)S(1,0)\cdots S(N,0)\cdots \le (d_{n+i})
\end{equation*}
by Lemma \ref{l32}.
\end{proof}

\begin{definition}
We say that an admissible sequence $d$ is \emph{of finite type} if $d=S_h(N,1)^{\infty}$ with some nonnegative integer $N$ and a finite sequence $h=(h_1,\ldots,h_N)$ of positive integers. Otherwise (including the case $d=0^{\infty}$) it is said to be \emph{of infinite type}.
\end{definition}

\begin{lemma}\label{l34}
Let $d=(d_i)\ne 1^{\infty}$ be an admissible sequence.\mbox{}
\smallskip

(a) If $(d_i)=S(N,1)^{\infty}$ is of finite type (then $N\ge 1$ because $d\ne 1^{\infty}$) and  $(d_i')=(d_{i+1+\ell_N-\ell_{N-1}})$, then
\begin{equation*}
(d_{n+i}')\ge (d_i')>(d_{1+i})
\end{equation*}
whenever $d_n'= 0$. Moreover, we have
\begin{align}
&(d_i')=S(1,0)\cdots S(N-1,0)S(N,1)^{\infty}\label{38}
\intertext{and}
&(d_{1+i})=S(1,0)\cdots S(N-1,0)S(N,0)S(N,1)^{\infty}.\label{39}
\end{align}
\smallskip

(b) In the other cases the sequence $(d_i'):=(d_{1+i})$ satisfies
\begin{equation*}
(d_{n+i}')\ge (d_i')
\end{equation*}
whenever $d_n'= 0$.
\smallskip

(c) We have $d'=d$ if and only if $d=\left( 1^{k-1}0\right)^{\infty} $ for some positive integer $k$, i.e., $d=0^{\infty}$ or $d=S(N,1)^{\infty}$ with $N=1$.
\end{lemma}

\begin{proof}\mbox{}

(a) The first inequality follows from Lemma \ref{l32}; the proof also shows that we have equality if and only if $n$ is a multiple of $\ell_N$.

The relations \eqref{32} and \eqref{33} of Lemma \ref{l31} imply \eqref{38}--\eqref{39} and they imply
the second inequality because $S(N,0)<S(N,1)$.
\smallskip

(b) The case $(d_i)=0^{\infty}$ is obvious. Otherwise $(d_i)$ begins with $S(N,1)$ for all $N\ge 0$ and $\ell_N\to\infty$, so that we deduce from the relation \eqref{32} of Lemma \ref{l31} the equality
\begin{equation*}
0d_2d_3\cdots =S(0,0)S(1,0)\cdots .
\end{equation*}
On the other hand, it follows from Lemma \ref{l32} that for any $n\ge 0$ we have
\begin{equation*}
(d_{n+i}')\ge S(0,0)\cdots S(N-1,0)S(N,0)^{\infty}
\end{equation*}
for every $N\ge 0$. This implies that
\begin{equation*}
(d_{n+i}')\ge 0d_2d_3\cdots
\end{equation*}
for every $n\ge 0$. If $d_n'= 0$, then we conclude that
\begin{equation*}
d'_nd'_{n+1}d'_{n+2}\cdots \ge 0d_2d_3\cdots
\end{equation*}
which is equivalent to the required inequality
\begin{equation*}
d'_{n+1}d'_{n+2}\cdots \ge d_2d_3\cdots .
\end{equation*}
\smallskip

(c) It follows from the above proof that $d=d'$ if and only if $d=0^{\infty}$ or  $d=S(N,1)^{\infty}$ for some integer $N\ge 1$ and $h$ such that $\ell_{N-1}=1$. These conditions are equivalent to  $d=\left( 1^{k-1}0\right)^{\infty} $ for some positive integer $k$.
\end{proof}

\begin{example}\label{e32}
If $d=S(N,1)^{\infty}$ is an admissible sequence of finite type with $N\ge 1$ (i.e., $d\ne 1^{\infty}$) and $h_1,\ldots, h_N\ge 1$, then there exists a smallest admissible sequence
$\tilde d>d$. It is of infinite type,  corresponding to the infinite sequence $h=h_1\cdots h_{N-1}h_N^+1^{\infty}$ with $h_N^+:=1+h_N$. Observe that $\tilde d=S(N-1,1)d$ and hence $\tilde d'=d'$.
\end{example}

\begin{lemma}\label{l35}
If $d=(d_i)$ is an admissible sequence, then no sequence $(c_i)$ of zeroes and ones satisfies
\begin{equation}\label{310}
0d_2d_3\cdots < (c_{n+i})< d_1d_2d_3\cdots
\end{equation}
for all $n=1,2,\ldots .$
\end{lemma}

\begin{proof}
The case $d=0^{\infty}$ is obvious because then $0d_2d_3\cdots = d_1d_2d_3\cdots .$ We may therefore assume that $d$ begins with $S(N,1)$ for some $N\ge 0$ and for some $h=(h_i)$. It follows from the second inequality of \eqref{310} that $(c_i)$ contains infinitely many zero digits. By removing a finite initial block from $(c_i)$ if necessary we may assume henceforth that $c_1=0$. Only these assumptions will be used in the first three steps below.
\medskip

\emph{First step.} The sequence $(c_i)$ cannot end with $S(N,0)^{\infty}$.

This is true for $h_1=\cdots =h_N=1$ because then we have $S(N,0)=0$ and obviously $0^{\infty}\le 0d_2d_3\cdots .$

Otherwise there exists $1\le M\le N$ such that $h_M\ge 2$ and $h_{M+1}=\cdots h_N=1$. Then we have $S(M,0)=S(M+1,0)=\cdots =S(N,0)$ and using the relation \eqref{34} of Lemma \ref{l31} there exists $n$ such that
\begin{align*}
(c_{n+i})
&=S(0,0)\cdots S(M-1,0)S(M,0)^{\infty}\\
&=0d_2d_3\cdots d_{\ell_N}S(N,0)^{\infty}\\
&<0d_2d_3\cdots
\end{align*}
because
\begin{equation*}
1d_2d_3\cdots d_{\ell_N}S(N,0)^{\infty}=S(N,1)S(N,0)^{\infty}
\end{equation*}
while $d$ begins with either with $S(N,1)S(N,1)$ or with $S(N,1)S(N,0)S(N,1)$ by Lemma \ref{l33}.
\medskip

\emph{Second step.}  We have $(c_i)=B_0\cdots B_N\left( S(N,c^N_j)\right) $ with some nonempty blocks $B_j\to S(j,0)$ and a suitable sequence $(c^N_j)\subset\set{0,1}$.

The case $N=0$ is obvious because $c_i=S(0,c_i)$ for every $i$; since $c_1=0$ by assumption we may choose $B_0=0$. Let $N\ge 1$ and assume by induction that $(c_i)=B_0\cdots B_{N-1}\left( S(N-1,c^{N-1}_j)\right) $ for some  blocks $B_j\to S(j,0)$ and a suitable sequence $(c^{N-1}_j)\subset\set{0,1}$.

Since $1d_2d_3\cdots$ begins with $S(N,1)=S(N-1,1)^{h_N}S(N-1,0)$, each block $S(N-1,c^{N-1}_j)$ is followed by at most $h_N$ consecutive blocks $S(N-1,1)$.

If $N=1$, then the above means that each digit is followed by at most $h_1$ consecutive one digits. On the other hand, since
$0d_2d_3\cdots$ begins with $01^{h_1-1}$, each zero digit is followed by at least $h_1-1$ consecutive one digits. This implies that $(c_i)=B_1\left( S(1,c^1_j)\right) $ for some  block $B_1\to S(1,0)$ and a suitable sequence $(c^1_j)\subset\set{0,1}$.

If $N\ge 2$, then, since $0d_2d_3\cdots$ begins with
\begin{multline*}
S(0,0)\cdots S(N-1,0)S(N,0)\\
=\left[ S(0,0)\cdots S(N-2,0)\right] S(N-1,0)S(N-1,1)^{h_N-1}S(N-1,0)
\end{multline*}
and since
\begin{equation*}
S(0,0)\cdots S(N-2,0)\to S(N-1,1)
\end{equation*}
by Lemma \ref{l31}, each block $S(N-1,1)S(N-1,0)$ in $\left( S(N-1,c^{N-1}_j)\right) $ is followed by at least $h_N-1$ consecutive blocks $S(N-1,1)$.

Since $(c_i)$ cannot end with $S(N-1,0)^{\infty}$  by the first step, we conclude that $\left( S(N-1,c^{N-1}_j)\right)=B_N\left( S(N,c^N_j)\right) $ for some  block $B_N\to S(N,0)$ and a suitable sequence $(c^N_j)\subset\set{0,1}$.
\medskip

\emph{Third step.}  We have
\begin{equation*}
(c_{n+i})=S(0,0)\cdots S(N-1,0)S(N,0)S(N,a_1)S(N,a_2)\cdots
\end{equation*}
for some $n$ and for a suitable  sequence $(a_j)\subset\set{0,1}$.

We already know from the second step that $(c_i)$ ends with $\left( S(N,c^N_j)\right) $ for a suitable sequence $(c^N_j)\subset\set{0,1}$, and that there are infinitely many blocks $S(N,1)$ by the first step. Since $(c_{n+i})<d\le S(N,1)^{\infty}$ for every $n\ge 0$, there are also infinitely many blocks $S(N,0)$. Therefore $(c_i)$ ends with
\begin{equation*}
S(N,1)S(N,0)S(N,a_1)S(N,a_2)\cdots
\end{equation*}
and our claim follows because $S(0,0)\cdots S(N-1,0)\to S(N,1)$ by Lemma \ref{l31}.
\medskip

\emph{Fourth step.}  We complete the proof of the lemma in the case $d=S(N,1)^{\infty}$. We deduce from the third step that
\begin{align*}
(c_{n+i})
&=S(0,0)\cdots S(N-1,0)S(N,0)S(N,a_1)S(N,a_2)\cdots\\
&\le S(0,0)\cdots S(N-1,0)S(N,0)S(N,1)^{\infty}\\
&=0d_2d_3\cdots ,
\end{align*}
contradicting our assumption on  $(c_i)$.
\medskip

\emph{Fifth step.}  We complete the proof of the lemma in the case where $d$ begins with $S(N,1)$ for every $N\ge 0$ for some $h=(h_i)$. We know from the second step that
$(c_i)=B_0\cdots B_N\left( S(N,c^N_j)\right) $ with some nonempty blocks $B_j\to S(j,0)$. If $B_j\ne S(j,0)$ for at least one $j$, then we conclude by using Lemma \ref{l31} that
\begin{equation*}
(c_i)<S(0,0)\cdots S(N,0)S(N,a_1)S(N,a_2)\cdots
\end{equation*}
for every sequence $(a_i)$. Since $0d_2d_3\cdots=S(0,0)\cdots S(N,0)S(N,a_1)S(N,a_2)\cdots $ for a suitable sequence  $(a_i)$, this contradicts the first inequality in \eqref{310}.

If $B_j=S(j,0)$ for all $j$, then $(c_i)$ begins  with  $S(0,0)\cdots S(N,0)$ for every $N\ge 0$. This implies that $(c_i)=0d_2d_3\cdots$, contradicting the first inequality in \eqref{310} again.
\end{proof}

\begin{lemma}\label{l36}
If $d=(d_i)\ne 1^{\infty}$ is an admissible sequence, then no sequence $(c_i)$ of zeroes and ones satisfies
\begin{equation}\label{311}
0(d_i')<(c_{n+i})<1(d_i')
\end{equation}
for all $n=1,2,\ldots .$
\end{lemma}

\begin{proof}
If $d=0^{\infty}$, then $d'=0^{\infty}$ and the result is obvious. It is sufficient therefore to consider the case where $d=S_h(N,1)^{\infty}$ for some $N\ge 1$ and for some $h=(h_i)$, for otherwise we have $d=1d'$ and we may apply the preceding lemma. It follows again from our assumptions that $(c_i)$ contains infinitely many zero digits, and we may assume that $c_1=0$.

Observe that putting $h^+:=(h_1,\ldots, h_{N-1},1+h_N)$, $1d'$ begins with $S_{h^+}(N,1)$. Therefore, repeating the first three steps of the preceding proof we obtain that
\begin{equation*}
(c_{n+i})=S_{h^+}(0,0)\cdots S_{h^+}(N-1,0)S_{h^+}(N,1)S_{h^+}(N,a_1)S_{h^+}(N,a_2)\cdots
\end{equation*}
for some $n$ and for a suitable  sequence $(a_j)\subset\set{0,1}$. This can be rewritten in the form
\begin{equation}\label{312}
(c_{n+i})=S_h(0,0)\cdots S_h(N-1,0)S_{h^+}(N,1)S_{h^+}(N,a_1)S_{h^+}(N,a_2)\cdots .
\end{equation}
Since
\begin{equation*}
1d'= S_{h^+}(N,1)S_h(N,1)^{\infty},
\end{equation*}
it follows from \eqref{311}--\eqref{312} that 
\begin{equation*}
(c_{n+i})=S_h(0,0)\cdots S_h(N-1,0)S_{h^+}(N,1)S_h(N,1)^{\infty}.
\end{equation*}
Since
\begin{equation*}
S_h(0,0)\cdots S_h(N-1,0)\to S_h(N,1),
\end{equation*}
we conclude that
\begin{equation*}
(c_{n'+i})=S_h(0,0)\cdots S_h(N-1,0)S_h(N,1)^{\infty}
\end{equation*}
for some $n'\ge 1$. This, however, contradicts \eqref{311} because the right side is equal to $0d'$.
\end{proof}

\section{$m$-admissible sequences}
\label{s4}

Throughout this section we fix an admissible sequence $d=(d_i)\ne 1^{\infty}$ and we define the sequence $d'=(d_i')$ as in Lemma \ref{l34}. Furthermore, for any given real number $m>1$ we denote by $\delta=(\delta_i)$ and $\delta'=(\delta_i')$ the sequences obtained from $d$ and $d'$ by the substitutions $1\to m$ and $0\to 1$. We define the numbers $p_m',p_m''>1$ by the equations
\begin{align}
\sum_{i=1}^{\infty}\frac{\delta_i}{(p_m')^i} =m-1\label{41}
\intertext{and}
\sum_{i=1}^{\infty}\frac{m-\delta_i'}{(p_m'')^i} =1\label{42}
\end{align}
and we put $p_m:=\max\set{p_m',p_m''}$.

Introducing the conjugate of $\delta$ by the formula $\overline{\delta'_i}:=m-\delta'_i$ we may also write \eqref{41} and \eqref{42} in the more economical form
\begin{equation*}
\pi_{p_m'}(\delta)=m-1
\quad\text{and}\quad
\pi_{p_m''}\left( \overline{\delta'}\right) =1.
\end{equation*}

Let us also introduce the number
\begin{equation*}
P_m:=1+\sqrt{\frac{m}{m-1}}.
\end{equation*}
A direct computation shows that $P_m>1$ can also be defined by any of the following equivalent conditions:

\begin{align}
&(P_m-1)^2=\frac{m}{m-1};\label{43}\\
&\frac{m}{P_m}+\frac{1}{P_m}\left( \frac{m}{P_m-1}-1\right) = m-1;\label{44} \\
&(m-1)P_m-m= \frac{m}{P_m-1}-1.\label{45}
\end{align}

We begin by investigating the dependence of $P_m$, $p_m'$ and $p_m''$ on $m$. The following two lemmas establish in particular part (b) of Theorem \ref{t11}.

\begin{lemma}\label{l41}\mbox{}

(a) The function $m\mapsto P_m$ is continuous and strictly decreasing in $(1,\infty)$.
\smallskip

(b) The function $m\mapsto p_m'-P_m$ is continuous and strictly decreasing in $(1,\infty)$, and it has a unique zero $m_d$.
\smallskip

(c) The function $m\mapsto p_m''-P_m$ is continuous and strictly increasing in $(1,\infty)$, and it has a unique zero $M_d$.
\smallskip

(d) The function $m\mapsto p_m'-p_m''$ is continuous and strictly decreasing in $(1,\infty)$, and it has a unique zero $\mu_d$.
\smallskip

(e) The function $m\mapsto p_m$ is continuous in $(1,\infty)$, strictly decreasing in $(1,\mu_d]$ and strictly increasing in $[\mu_d,\infty)$, so that it has a strict global minimum in $\mu_d$.
\end{lemma}

\begin{proof}\mbox{}

(a) A straightforward computation shows that $P$ is infinitely differentiable in $(1,\infty)$ and
\begin{equation*}
P'(m)=\frac{-1}{2(m-1)\sqrt{m(m-1)}}<0
\end{equation*}
for all $m>1$.
\smallskip

(b) Since $\delta_i=1+(m-1)d_i$, we may rewrite \eqref{41} in the form
\begin{equation}\label{46}
\frac{1}{m-1}+(p_m'-1) \sum_{i=1}^{\infty} \frac{d_i}{(p_m')^i} =p_m'-1.
\end{equation}
Applying the implicit function theorem it follows that the function $m\mapsto p_m'$ is $C^{\infty}$.

Differentiating the last identity with respect to $m$, denoting the derivatives by dots and setting
\begin{equation*}
A:=1+(p_m'-1)\left( \sum_{i=1}^{\infty}\frac{d_i}{i(p_m')^{i+1}}\right)
-\left( \sum_{i=1}^{\infty} \frac{d_i}{(p_m')^i} \right),
\end{equation*}
we get
\begin{equation*}
A\dot p_m'=\frac{-1}{(m-1)^2}.
\end{equation*}
Differentiating \eqref{43} we obtain that the right side is equal to $2(P_m-1)\dot P_m$, so that
\begin{equation*}
A\dot p_m'=2(P_m-1)\dot P_m.
\end{equation*}
Since $\dot P_m<0$ and $2(P_m-1)>1$, it suffices to show that $A\in (0,1)$. Indeed, then we will have $\dot p_m'/\dot P_m>1$ and therefore $\dot p_m'<\dot P_m (<0)$.

The inequality $A>0$ follows by using \eqref{46}:
\begin{equation*}
A=(p_m'-1)\left( \sum_{i=1}^{\infty}\frac{d_i}{i(p_m')^{i+1}}\right)+\frac{1}{(m-1)(p_m'-1)}>0,
\end{equation*}
while the proof of $A<1$ is straightforward:
\begin{align*}
A
&\le 1+(p_m'-1)\left( \sum_{i=1}^{\infty}\frac{d_i}{(p_m')^{i+1}}\right)-\left( \sum_{i=1}^{\infty} \frac{d_i}{(p_m')^i} \right)\\
&=1-\frac{1}{p_m'}\left( \sum_{i=1}^{\infty} \frac{d_i}{(p_m')^i} \right)\\
&<1.
\end{align*}

It remains to show that $p_m'-P_m$ changes sign in $(1,\infty)$. It is clear from the definition that
\begin{equation}\label{47}
\lim_{m\searrow 1}P_m=\infty
\quad\text{and}\quad
\lim_{m\to\infty}P_m=2.
\end{equation}
Furthermore, using the equality $d_1=1$ it follows from \eqref{46} that
\begin{equation*}
\frac{1}{m-1}\le p_m'-1\le 1+\frac{1}{m-1};
\end{equation*}
hence
\begin{equation}\label{48}
\lim_{m\searrow 1}p_m'=\infty
\quad\text{and}\quad
\lim_{m\to\infty}p_m'=1.
\end{equation}
We infer from \eqref{47}--\eqref{48} that $\lim_{m\to\infty}p_m'-P_m=-1<0$. The proof is completed by observing that
\begin{equation*}
p_m'-P_m\ge \frac{1}{m-1}-1-\sqrt{\frac{m}{m-1}}\to \infty>0
\end{equation*}
if $m\searrow 1$.
\smallskip

(c) We may rewrite \eqref{42} in the form
\begin{equation}\label{49}
\sum_{i=1}^{\infty} \frac{1-d_i'}{(p_m'')^i} =\frac{1}{m-1}.
\end{equation}
Applying the implicit function theorem it follows from \eqref{49} that the function $m\mapsto p_m''$ is $C^{\infty}$.

The last identity also shows that the function $m\mapsto p_m''$ is strictly increasing. Using (a) we conclude that the function $m\mapsto p_m''-P_m$ is strictly increasing, too.

It remains to show that $p_m''-P_m$ changes sign in $(1,\infty)$. Since $d\ne 1^{\infty}$, there exists an index $k$ such that $d_k'=0$. Therefore we deduce from \eqref{49} the inequalities
\begin{equation*}
\frac{1}{(p_m'')^k}\le \frac{1}{m-1}\le \frac{1}{p_m''-1}
\end{equation*}
and hence that
\begin{equation}\label{410}
\lim_{m\searrow 1}p_m''=1
\quad\text{and}\quad
\lim_{m\to\infty}p_m''=\infty.
\end{equation}
We conclude from \eqref{47} and \eqref{410} that
\begin{equation*}
\lim_{m\searrow 1}p_m''-P_m=-\infty<0
\quad\text{and}\quad
\lim_{m\to\infty}p_m''-P_m=\infty>0.
\end{equation*}
\smallskip

(d) The proof of (b) and (c) shows that $m\mapsto p_m'$ is continuous and strictly decreasing and $m\mapsto p_m''$ is continuous and strictly increasing; hence the function $m\mapsto p_m'-p_m''$ is continuous and strictly decreasing.  It remains to observe that $p_m'-p_m''$ changes sign in $(1,\infty)$ because \eqref{48} and \eqref{410} imply that
\begin{equation*}
\lim_{m\searrow 1}p_m'-p_m''=\infty>0
\quad\text{and}\quad
\lim_{m\to\infty}p_m'-p_m''=-\infty<0.
\end{equation*}
\smallskip

(e) This follows from the definition $p_m:=\max\set{p_m',p_m''}$ and from the fact that $m\mapsto p_m'$ is continuous and strictly decreasing and $m\mapsto p_m''$ is continuous and strictly increasing.
\end{proof}

The following lemma is a variant of a similar result in \cite{EJK1}.

\begin{lemma}\label{l42}
Let $(c_i)$ be an expansion of some number $s\le b-a$ in some base $q>1$ on some alphabet $\set{a,b}$ with $a<b$. If 
\begin{equation*}
c_{n+1}c_{n+2}\cdots\le c_1c_2\cdots\quad\text{whenever}\quad c_n=a,
\end{equation*}
then we also have
\begin{equation*}
\frac{c_{n+1}}{q^{n+1}}+\frac{c_{n+2}}{q^{n+2}}+\cdots \le \frac{s}{q^n}\quad\text{whenever}\quad c_n=a.
\end{equation*}

Moreover, the inequality is strict if the sequence $(c_i)$ is infinite and $(c_{n+i})\ne (c_i)$.
\end{lemma}

\begin{proof}
Starting with $k_0:=n$ we define by recurrence a sequence of indices
$k_0<k_1<\cdots$ satisfying for $j=1,2,\ldots$ the conditions
\begin{equation*}
c_{k_{j-1}+i}=c_i\quad\text{for}\quad i=1,\ldots, k_j-k_{j-1}-1,
\quad\text{and}\quad c_{k_j}<c_{k_j-k_{j-1}}.
\end{equation*}
If we obtain an infinite sequence, then we have
\begin{align*}
\sum _{i=n+1}^\infty \frac{c_i}{q^i}
&= \sum _{j=1}^\infty\sum_{i=1}^{k_j-k_{j-1}}\frac{c_{k_{j-1}+i}}{q^{k_{j-1}+i}}\\
&\le \sum _{j=1}^\infty \Bigl(\Bigl(\sum_{i=1}^{k_j-k_{j-1}}\frac{c_i}{q^{k_{j-1}+i}}\Bigr)-\frac{b-a}{q^{k_j}}\Bigr)\\
&\le\sum _{j=1}^\infty \Bigl(\frac{s}{q^{k_{j-1}}}-\frac{b-a}{q^{k_j}}\Bigr)\notag\\
&\le\sum _{j=1}^\infty \Bigl(\frac{s}{q^{k_{j-1}}}-\frac{s}{q^{k_j}}\Bigr)\notag\\
&=\frac{s}{q^n}.
\end{align*}

Otherwise we have $(c_{k_N+i})=(c_i)$ after a finite number of steps (we do
not exclude the possibility that $N=0$), and we may conclude as follows:
\begin{align*}
\sum _{i=n+1}^\infty \frac{c_i}{q^i}
&=\Bigl(\sum _{j=1}^N\sum_{i=1}^{k_j-k_{j-1}}\frac{c_i}{q^{k_{j-1}+i}}\Bigr)+\sum _{i=1}^\infty \frac{c_{k_N+i}}{q^{k_N+i}}\\
&\le \sum _{j=1}^N \Bigl(\Bigl(\sum_{i=1}^{k_j-k_{j-1}}\frac{c_i}{q^{k_{j-1}+i}}\Bigr)-\frac{b-a}{q^{k_j}}\Bigr)+\sum_{i=1}^\infty\frac{c_i}{q^{k_N+i}}\\
&\le\sum _{j=1}^N \Bigl(\frac{s}{q^{k_{j-1}}}-\frac{b-a}{q^{k_j}}\Bigr)+\frac{s}{q^{k_N}}\\
&\le\sum _{j=1}^N \Bigl(\frac{s}{q^{k_{j-1}}}-\frac{s}{q^{k_j}}\Bigr)+\frac{s}{q^{k_N}}\\
&=\frac{s}{q^n}.
\end{align*}

The last property follows from the above proof.
\end{proof}

Now we investigate the mutual positions of $m_d$, $M_d$ and $\mu_d$.

\begin{lemma}\label{l43}\mbox{}

(a) If $d$ is of finite type, then $m_d<\mu_d<M_d$, and $p_m< P_m$ for all $m_d<m<M_d$. Furthermore, $p_m\ge 2$ for all $m\in (1,\infty)$ with equality if and only if $d=\left( 1^{k-1}0\right)^{\infty} $ and $m=2^k$ for some positive integer $k$.
\smallskip

(b) In the other cases we have $m_d=\mu_d=M_d$ and $p_m\ge p_{\mu_d}=P_{\mu_d}>2$ for all $m\in (1,\infty)$.
\end{lemma}

\begin{proof}\mbox{}

(a) In view of Lemma \ref{l41} the first assertion will follow if we show that $p_m< P_m$ for $m:=\mu_d$. Since $p_m=p_m'=p_m''$ in this case, using the relations \eqref{38}--\eqref{39} of Lemma \ref{l34} we have
\begin{align*}
m-1
&= \sum_{i=1}^{\infty} \frac{\delta_i}{p_m^i} \\
&=\frac{m}{p_m}+\frac{1}{p_m}\sum_{i=1}^{\infty} \frac{\delta_{i+1}}{p_m^i} \\
&< \frac{m}{p_m}+\frac{1}{p_m}\sum_{i=1}^{\infty} \frac{\delta_i'}{p_m^i} \\
&=\frac{m}{p_m}+\frac{1}{p_m}\left( \frac{m}{p_m-1}-1\right).
\end{align*}
In this computation the crucial inequality follows from Lemmas \ref{l34} and \ref{l42}. Indeed, writing $d=S(N,1)^{\infty}$, in view of the relations \eqref{38}--\eqref{39} of Lemma \ref{l34} the inequality is equivalent to
\begin{equation*}
\pi_{p_m'}\left( (\delta_{\ell_{N-1}+i})\right) <\pi_{p_m'}(\delta),
\end{equation*}
and this inequality follows from Lemma \ref{l42} with $c=\delta$, $q=p_m'$ and $n=\ell_{N-1}$. (The hypotheses of the lemma are satisfied because $d$ is an admissible sequence.)

Using \eqref{44} we conclude that $p_m< P_m$ indeed.

Furthermore, for $m:=\mu_d$ we deduce from the equalities
\begin{equation*}
\pi_{p_m}(\delta)=m-1
\quad\text{and}\quad
\pi_{p_m}\left( \overline{\delta'}\right) =1
\end{equation*}
that
\begin{equation*}
\sum_{i=1}^{\infty} \frac{m-\delta'_i+\delta_i}{p_m^i} =m.
\end{equation*}
It follows that $p_m\ge 2$ if and only if
\begin{equation*}
\sum_{i=1}^{\infty} \frac{m-\delta'_i+\delta_i}{2^i}\ge m
\end{equation*}
which is equivalent to the inequality
\begin{equation*}
\pi_2(\delta')\le \pi_2(\delta).
\end{equation*}
Since $\delta'\le \delta$ by Lemma \ref{l34}, this is satisfied by a well-known property of diadic expansions.

The proof also shows that we have equality if and only if $\delta'=\delta$. By part (c) of Lemma \ref{l34} this is equivalent to  $d=\left( 1^{k-1}0\right)^{\infty} $ for some positive integer $k$. In this case we infer from the equations
\begin{align*}
&\frac{m}{p_m'-1}-\frac{m-1}{(p_m')^k-1}=m-1
\intertext{and}
&\frac{m}{p_m''-1}-\frac{m-1}{(p_m'')^k-1}=\frac{m}{p_m''-1}-1
\end{align*}
that $p_m'=p_m''=m^{1/k}=2$.

Since by Lemma \ref{l41} $p_m$ has a global strict minimum in $m=\mu_d$, we have $p_m>2$ for all other values of $m$.
\smallskip

(b) Putting $m=\mu_d$ and repeating the first part of the proof of (a), by Lemma \ref{l34} now we have an equality instead of the strict inequality; using \eqref{44} we conclude that $p_m=P_m$ and hence $p_m=p_m'=p_m''=P_m$. Applying Lemma \ref{l41} we conclude that $m_d=\mu_d=M_d$.
\end{proof}

\section{Univoque sequences in small bases}
\label{s5}

In this section we determine the generalized golden ratio for every ternary alphabet $A=\set{a_1,a_2,a_3}$. Putting
\begin{equation*}
m:=\max\left\lbrace \frac{a_3-a_1}{a_2-a_1},\frac{a_3-a_1}{a_3-a_2}\right\rbrace
\end{equation*}
we will show that
\begin{equation*}
2\le G_A\le P_m:=1+\sqrt{\frac{m}{m-1}}.
\end{equation*}
Moreover, we will give an exact expression of $G_A$ for each $m$ and we will determine the values of $m$ for which $G_A=2$ or $G_A=P_m$.

By Lemma \ref{l25} we may restrict ourselves without loss of generality to the case of the alphabets $A_m=\set{0,1,m}$ with $m\ge 2$. Condition \eqref{21} takes the form
\begin{equation*}
1<q\le \frac{2m-1}{m-1};
\end{equation*}
under this assumption, that we assume henceforth, the results of the preceding section apply. For the sequel we fix a real number $m\ge 2$ and we consider expansions in bases $q>1$ with respect to the ternary alphabet $A_m:=\set{0,1,m}$ .

One of our main tools will be Theorem \ref{t22} which now takes the following special form:

\begin{lemma}\label{l51}
An expansion $(c_i)$ is unique in base $q$ for the alphabet $A_m$ if and only if the following conditions are satisfied:
\begin{align}
&\sum_{i=1}^\infty \frac{c_{n+i}}{q^i}  < 1 \quad\text{whenever $c_n=0$};\label{51}\\
&\sum_{i=1}^\infty \frac{c_{n+i}}{q^i}  < m-1 \quad\text{whenever $c_n=1$};\label{52}\\
&\sum_{i=1}^\infty \frac{c_{n+i}}{q^i}  > \frac{m}{q-1}-1 \quad\text{whenever $c_n=1$};\label{53}\\
&\sum_{i=1}^\infty \frac{c_{n+i}}{q^i}  > \frac{m}{q-1}-(m-1) \quad\text{whenever $c_n=m$}.\label{54}
\end{align}
\end{lemma}

\begin{corollary}\label{c52}
We have $G_{A_m}\ge 2$.
\end{corollary}

\begin{proof}
Let $(c_i)$ be a univoque sequence in some base $1<q\le 2$. We infer from \eqref{52} and \eqref{53} that $c_n\ne 1$ for every $n$. Since $m\ge q$, then we conclude from \eqref{51} that each $0$ digit is followed by another $0$ digit. Therefore condition \eqref{54} implies that each $m$ digit is followed by another $m$ digit. For otherwise the left-hand side of \eqref{54} would be zero, while the right-hand side is equal to one. Hence $(c_i)$ must be equal to $0^{\infty}$ or $m^{\infty}$.
\end{proof}

\begin{lemma}\label{l53}
If $(c_i)$ is a nontrivial univoque sequence in some base $1<q\le P_m$, then $(c_i)$ contains at most finitely many zero digits.
\end{lemma}

\begin{proof}
Since a univoque sequence remains univoque in every larger base, too, we may assume that $q=P_m$.
It suffices to prove that $(c_i)$ does not contain any block of the form $m0$ or $10$.

\emph{$(c_i)$ does not contain any block of the form $m0$.}  If $c_n=m$ and $c_{n+1}=0$ for some $n$, then we deduce from Lemma \ref{l51} that
\begin{equation*}
\sum_{i=1}^\infty \frac{c_{n+i}}{P_m^i} >\frac{m}{P_m-1}-(m-1)
\quad\text{and}\quad
\sum_{i=1}^\infty \frac{c_{n+i+1}}{P_m^i}  < 1.
\end{equation*}
Hence
\begin{equation*}
\frac{m}{P_m-1}-(m-1)<\sum_{i=1}^\infty \frac{c_{n+i}}{P_m^i} =\frac{1}{P_m}\sum_{i=1}^\infty \frac{c_{n+i+1}}{P_m^i}<\frac{1}{P_m},
\end{equation*}
contradicting one of our conditions on $P_m$ above.
\smallskip

\emph{$(c_i)$ does not contain any block of the form $10$.} If $c_n=1$ and $c_{n+1}=0$ for some $n$, then the application of Lemma \ref{l51} shows that
\begin{equation*}
\sum_{i=1}^\infty \frac{c_{n+i}}{P_m^i} >\frac{m}{P_m-1}-1
\quad\text{and}\quad
\sum_{i=1}^\infty \frac{c_{n+i+1}}{P_m^i}  < 1.
\end{equation*}
Since $m\ge 2$, these inequalities imply those of the preceding step, contradicting again our condition on $P_m$.
\end{proof}

Next we select a particular admissible sequence for each given $m$. Given an admissible sequence $d\ne 1^{\infty}$ we set
\begin{equation}\label{55}
I_d:=
\begin{cases}
[m_d,M_d)&\text{if $m_d<M_d$},\\
\set{m_d}&\text{if $m_d=M_d$}.
\end{cases}
\end{equation}

\begin{lemma}\label{l54}
Given a real number $m\ge 2$ there exists a lexicographically largest admissible sequence $d=(d_i)$ such that using the notations of the preceding section we have
\begin{equation}\label{56}
\sum_{i=1}^\infty \frac{\delta_i}{P_m^i} \le m-1.
\end{equation}

Furthermore, we have $d\ne 1^{\infty}$ and $m\in I_d$.
\end{lemma}

\begin{remark}
The lemma and its proof remain valid for all $m\ge (1+\sqrt{5})/2$.
\end{remark}

\begin{proof}
The sequence $d=0^{\infty}$ always satisfies \eqref{56} because
\begin{equation*}
\frac{1}{P_m-1}\le 1
\Longleftrightarrow
\sqrt{\frac{m-1}{m}}\le m-1
\Longleftrightarrow
1\le m(m-1)
\end{equation*}
and the last inequality is satisfied for all $m\ge 2$ (and even for all $m\ge (1+\sqrt{5})/2$).

All other admissible sequences are defined by a finite or infinite sequence $(h_j)$. If we add the symbol $\infty$ to the end of each finite sequence  $(h_j)$, then the lexicographic order between admissible sequences is equivalent to the lexicographic order between the corresponding sequences $(h_j)$. Let us say that a sequence $(h_j)$ is \emph{suitable} if the corresponding admissible sequence satisfies \eqref{56}. We are thus looking for the largest suitable sequence  $(h_j)$.

If $1{\infty}$ is not suitable, then $d=0^{\infty}$ is the largest admissible sequence satisfying \eqref{56}.

If $1{\infty}$ is suitable, then there exists a largest positive integer $h_1$ for which $h_11^{\infty}$ is suitable. For otherwise $d=1^{\infty}$
would satisfy  \eqref{56} which is impossible because for $d=1^{\infty}$ the right side of  \eqref{56} is equal to
\begin{equation*}
\frac{m}{P_m-1}=m-1+\frac{1}{P_m}>m-1.
\end{equation*}

Proceeding by recurrence assume that we have already determined the largest positive integers $h_1$,\ldots, $h_N$ for some $N\ge 1$ such that $h_1\cdots h_N1^{\infty}$ is suitable. If $h_1\cdots h_N{\infty}$ is suitable, then it is the largest suitable sequence. If not, then there exists a largest positive integer $h_{N+1}$ such that $h_1\cdots h_Nh_{N+1}1^{\infty}$ is suitable.

Continuing this recurrence, either we find a largest suitable sequence of the form $h_1\cdots h_N{\infty}$ with $N\ge 1$ in a finite number of steps, or we construct a largest infinite suitable sequence  $(h_j)$.

Observe that $d\ne 1^{\infty}$ because for $d= 1^{\infty}$ using \eqref{43} we have
\begin{equation*}
\sum_{i=1}^\infty \frac{\delta_i}{P_m^i} =\frac{m}{P_m-1}=\sqrt{(m-1)m} > m-1
\end{equation*}
so that \eqref{56} is not satisfied.

It remains to prove that $m\in I_d$. We distinguish three cases.

(a) If $(d_i)$ is defined by an infinite sequence $(h_j)$, then we already know that $p_m=p_m'=p_m''$ and that
\begin{equation*}
\sum_{i=1}^\infty \frac{\delta_i}{P_m^i} \le m-1.
\end{equation*}
It remains to show the converse inequality
\begin{equation}\label{57}
\sum_{i=1}^\infty \frac{\delta_i}{P_m^i}\ge m-1.
\end{equation}
It follows from the definition of $(\delta_i)$ that if we denote by $(\delta^N_i)$ the sequence associated with the admissible sequence defined by the sequence $h:=h_1,\ldots,h_{N-1},h_N+1,1,1,\ldots$, then
\begin{equation*}
\sum_{i=1}^\infty \frac{\delta_i^N}{P_m^i} > m-1.
\end{equation*}
Since both $(d_i)$ and $(d^N_i)$ begin with $S(N-1,1)^{h_N}$ and since the length of this block tends to infinity, letting $N\to\infty$ we conclude \eqref{57}.
\smallskip

(b) If $(d_i)=S(N,1)^{\infty}$ for some $N\ge 1$, then
\begin{align*}
(e_i):&=S(N-1,1)^{h_N+1}S(N-1,0)\left[ S(N-1,1)^{h_N}S(N-1,0)\right] ^{\infty}\\
&=S(N-1,1)S(N,1)^{\infty}
\end{align*}
does not satisfy \eqref{56}, so that
\begin{equation*}
\sum_{i=1}^\infty \frac{\eps_i}{P_m^i} > m-1
\end{equation*}
where $(\eps_i)$ is obtained from $(e_i)$ by the usual substitutions $1\to m$ and $0\to 1$.

Observe that now we have $e_1e_2\cdots = 1d_1'd_2'\cdots$ and therefore (using the notations of the first page of the paper)
\begin{equation*}
m-1<\pi_{P_m}(\eps)=\frac{m}{P_m}+\frac{1}{P_m}\pi_{P_m}(\delta').
\end{equation*}
It follows that
\begin{equation*}
\pi_{P_m}(\delta')>(m-1)P_m-m=\frac{m}{P_m-1}-1
\end{equation*}
which is equivalent to $\pi_{P_m}\left( \overline{\delta'}\right) <1$. Since we have $\pi_{p_m''}\left( \overline{\delta'}\right) =1$ by the definition of $p_m''$, we conclude that $P_m>p_m''$.

Finally, since we have $\pi_{P_m}(\delta)\le m-1=\pi_{p_m'}(\delta)$ by the definitions of $(d_i)$ and $p_m'$, we have also  $P_m\ge p_m'$.
\smallskip

(c) If $(d_i)=0^{\infty}$, then we may repeat the proof of (b) with $(d_i')=0^{\infty}$ and $(e_i)=10^{\infty}$.
\end{proof}

\begin{example}\label{e51}
Using a computer program we can determine the admissible sequences of Lemma \ref{l54} for all \emph{integer} values $2\le m\le 2^{16}$. For all but seven values the corresponding admissible sequence is of finite type with $N=1$, more precisely $d=(1^{h_1}0)^{\infty}$ with $h_1=[\log_2 m]$. For the exceptional values $m=5,9,130,258,2051, 4099, 32772$ the corresponding admissible sequence is of finite type with $N=2$ and $h_1=[\log_2 m]$ as shown  in the following table:
\medskip

\begin{center}
\begin{tabular}{|c|c|c|c|}
\hline
$m$ & $d$ & $N$ & $h$\\
\hline
$5$ & $(11011010)^\infty$ & 2 & (2,2)\\
$9$ & $(1110110)^\infty$ & 2 & (3,1)\\
$130$ & $(1^701^60)^\infty$ &  2 & (7,1)\\
$258$ & $(1^801^70)^\infty$ &  2 & (8,1)\\
$2051$ & $(1^{11}01^{10}0)^\infty$ &  2 & (11,1)\\
$4099$ & $(1^{12}01^{11}0)^\infty$ &  2 & (12,1)\\
$32772$ & $(1^{15}01^{14}0)^\infty$ &  2 & (15,1)\\
\hline
\end{tabular}
\end{center}
\medskip

\end{example}

Now we need two definitions. The \emph{quasi-greedy} expansion of a real number $x$ in some base $q$ is its lexicographically largest infinite expansion in the alphabet $\set{0,1,m}$, while the \emph{quasi-lazy} expansion of $x$ is the \emph{conjugate} $(m-c_i)$ of the quasi-greedy expansion $(c_i)$ of $\frac{m}{q-1}-x$ with respect to the conjugate alphabet $\set{0,m-1,m}$. The following lemma follows at once from these definitions.

\begin{lemma}\label{l55}
Let $(c_i)$ be a sequence on the alphabet $\set{0,1,m}$ and $q>1$ a real number.\mbox{}

(a) The sequence $(c_i)$ is a quasi-greedy expansion of some $x$ in base $q$ if and only if
\begin{align*}
&\sum_{i=1}^{\infty}\frac{c_{n+i}}{q^i}\le 1\quad\text{whenever}\quad c_n=0
\intertext{and}
&\sum_{i=1}^{\infty}\frac{c_{n+i}}{q^i}\le m-1\quad\text{whenever}\quad c_n=1.
\end{align*}
Hence, if $c=(c_i)$ is a quasi-greedy expansion in base $q$, then $m^nc$ and $(c_{n+i})$ are also quasi-greedy expansions in every base $\ge q$, for every positive integer $n$.
\smallskip

(b) The sequence $(c_i)$ is a quasi-lazy expansion of some $x$ in base $q$ if and only if
\begin{align*}
&\sum_{i=1}^{\infty}\frac{m-c_{n+i}}{q^i}\le 1\quad\text{whenever}\quad c_n=1
\intertext{and}
&\sum_{i=1}^{\infty}\frac{m-c_{n+i}}{q^i}\le m-1\quad\text{whenever}\quad c_n=m.
\end{align*}
Hence, if $c=(c_i)$ is a quasi-lazy expansion in base $q$, then $0^nc$ and $(c_{n+i})$ are also quasi-lazy expansions in every base $\ge q$, for every positive integer $n$.
\smallskip

(c) If $x\ge y$ and $p\ge q$, then the quasi-greedy (resp. the quasi-lazy) expansion of $x$ in base $p$ is lexicographically larger than or equal to that of $y$ in base $q$.
\end{lemma}

\begin{lemma}\label{l56}
Given an admissible sequence $d\ne 1^{\infty}$ and $m\in I_d$ let us define the sequences $d'$, $\delta$, $\delta'$ and the numbers $p_m',p_m'',p_m$ as at the beginning of Section \ref{s4}.
\smallskip

(a) The sequences $\delta$, $m \delta'$ are quasi-greedy  in base $p_m$.
\smallskip

(b) The sequences $\delta'$ and $(\delta_{1+i})$ are quasi-lazy in base $p_m$.
\end{lemma}

\begin{proof}\mbox{}

(a) Applying Lemma \ref{l42} with $(c_i):=\delta$ and $q=p_m'$ on the alphabet $\set{1,m}$ we obtain that $\delta$ satisfies the conditions of part (a) Lemma \ref{l55}, so that it is a quasi-greedy expansion in base $p_m'$. Since $p_m\ge p_m'$, the second half of part (a) of Lemma \ref{l55} shows that $\delta$ and $m \delta'$ are quasi-greedy expansions in base $p_m$.
\smallskip

(b) Since $(\delta_{1+i})=(\delta'_{k+i})$ for some $k\ge 0$, and since $1\le m-1$, in view of part (b) of Lemma \ref{l55} it suffices to show that
\begin{equation}\label{58}
\sum_{i=1}^{\infty}\frac{m-\delta'_{n+i}}{(p_m'')^i}\le 1
\end{equation}
for all $n\ge 1$.

First we observe that \eqref{58} is true for $n=0$ by the definition of $p_m''$. Furthermore, it is also true if $n\ge 1$ and $\delta'_n=m$: this follows by applying Lemma \ref{l42} with $(c_i):=m-\delta'$ and $q=p_m''$ on the alphabet $\set{0,m-1}$.

If $n\ge 1$ and $\delta'_n=1$, then let $0\le k<n$ be the largest integer satisfying $\delta_i=1$ for all $k<i\le n$. Then we have
\begin{equation*}
\sum_{i=1}^{\infty}\frac{m-\delta'_{k+i}}{(p_m'')^i}\le 1
\end{equation*}
by the preceding paragraph. The proof will be complete if we show that
\begin{equation*}
\sum_{i=1}^{\infty}\frac{m-\delta'_{n+i}}{(p_m'')^i}\le \sum_{i=1}^{\infty}\frac{m-\delta'_{k+i}}{(p_m'')^i}.
\end{equation*}
Since
\begin{equation*}
\sum_{i=1}^{\infty}\frac{m-\delta'_{k+i}}{(p_m'')^i}=\frac{m-1}{p_m''}+\cdots+\frac{m-1}{(p_m'')^{n-k}}+\frac{1}{(p_m'')^{n-k}}\sum_{i=1}^{\infty}\frac{m-\delta'_{n+i}}{(p_m'')^i},
\end{equation*}
this is equivalent to
\begin{equation*}
\left( 1-\frac{1}{(p_m'')^{n-k}}\right) \sum_{i=1}^{\infty}\frac{m-\delta'_{n+i}}{(p_m'')^i}\le \frac{m-1}{p_m''}+\cdots+\frac{m-1}{(p_m'')^{n-k}}.
\end{equation*}
The last inequality is true because if we replace each $\delta'_{n+i}$ by the smallest possible value $1$, then we obtain the right side. Indeed, 
\begin{align*}
\left( 1-\frac{1}{(p_m'')^{n-k}}\right) \sum_{i=1}^{\infty}\frac{m-1}{(p_m'')^i}
&=\left(  \sum_{i=1}^{\infty}\frac{m-1}{(p_m'')^i}\right) -\left(  \sum_{i=n-k+1}^{\infty}\frac{m-1}{(p_m'')^i}\right)\\
&=\frac{m-1}{p_m''}+\cdots+\frac{m-1}{(p_m'')^{n-k}}.\qedhere
\end{align*}
\end{proof}

\begin{lemma}\label{l57}
Denoting by $\gam=(\gam_i)$ and $\lam=(\lam_i)$ the quasi-greedy expansion of $m-1$  in base $p_m$ and the quasi-lazy expansion of $\frac{m}{p_m-1}-1$ in base $p_m$, respectively, we have either
\begin{equation*}
(\delta_{1+i})\le \lam \quad\text{and}\quad \gam=\delta
\end{equation*}
or
\begin{equation*}
\delta'=\lam \quad\text{and}\quad \gam \le m\delta' .
\end{equation*}
\end{lemma}

\begin{proof}
If $p_m'\ge p_m''$, then both $\gam$ and $\delta $ are quasi-greedy expansions of $m-1$ in base $p_m=p_m'$ by Lemma \ref{l56}, so that $\gam =\delta $. Since furthermore both $\hat\delta:=(\delta_{1+i})$ and $\lam $ are quasi-lazy expansions in base $p_m$, in view of Lemma \ref{l55} it remains to show only that $\pi_{p_m}(\hat\delta)\le \pi_{p_m}(\lam )$. Since
\begin{equation*}
m-1=\pi_{p_m}(\delta)=\frac{m}{p_m}+\frac{1}{p_m}\pi_{p_m}(\hat\delta)
\end{equation*}
and $p_m\le P_m$, using \eqref{45} we have
\begin{equation*}
\pi_{p_m}(\hat\delta)=(m-1)p_m-m\le \frac{m}{p_m-1}-1=\pi_{p_m}(\lam ).
\end{equation*}

If  $p_m''\ge p_m'$, then both $\lam$ and $\delta' $ are quasi-lazy expansions of $\frac{m}{p_m-1}-1$ in base $p_m=p_m''$ by Lemma \ref{l56}, so that  $\lam =\delta'$. Furthermore $m\delta'$ and $\gam $ are quasi-greedy expansions in base $p_m$. Since $p_m\le P_m$, using \eqref{44} we obtain that
\begin{align*}
\pi_{p_m}(m\delta')
&=\frac{m}{p_m}+\frac{1}{p_m}\pi_{p_m}(\delta')\\
&=\frac{m}{p_m}+\frac{1}{p_m}\left( \frac{m}{p_m-1}-1\right) \\
&\ge m-1\\
&=\pi_{p_m}(\gam ).
\end{align*}
Applying Lemma \ref{l55} we conclude that $m\delta'\ge\gam $.
\end{proof}

Given $m\ge 2$ we choose an admissible sequence $d\ne 1^{\infty}$ satisfying $m\in I_d$ (see Lemma \ref{l54}) and we define $p_m$ as at the beginning of Section \ref{s4} (see Lemma \ref{l56}). The following lemma proves part (a) of Theorem \ref{t11}.

\begin{lemma}\label{l58}\mbox{}

(a) If $q>p_m$, then $\delta'$ is a nontrivial univoque sequence in base $q$.
\smallskip

(b) There are no nontrivial univoque sequences in any base $1<q\le p_m$.
\end{lemma}

\begin{proof}\mbox{}

(a) Since the sequence $\delta$ is quasi-greedy and the sequence $\delta'$ is quasi-lazy in base $p_m$ and since $\delta'$ is obtained from $\delta$ by removing a finite initial block, $\delta'$ is both quasi-greedy and quasi-lazy in base $p_m$. Hence
\begin{align*}
&\sum_{i=1}^\infty \frac{\delta_{n+i}'}{(p_m)^i} \le m-1\quad\text{whenever}\quad \delta_n'=1,\\
&\sum_{i=1}^\infty \frac{m-\delta_{n+i}'}{(p_m)^i}  \le 1\quad\text{whenever}\quad \delta_n'=1,\\
&\sum_{i=1}^\infty \frac{m-\delta_{n+i}'}{(p_m)^i}  \le m-1\quad\text{whenever}\quad \delta_n'=m.
\end{align*}
Since $q>p_m$, it follows that
\begin{align*}
&\sum_{i=1}^\infty \frac{\delta_{n+i}'}{q^i} < m-1\quad\text{whenever}\quad \delta_n'=1,\\
&\sum_{i=1}^\infty \frac{m-\delta_{n+i}'}{q^i}  < 1\quad\text{whenever}\quad \delta_n'=1,\\
&\sum_{i=1}^\infty \frac{m-\delta_{n+i}'}{q^i} < m-1\quad\text{whenever}\quad \delta_n'=m.
\end{align*}
Applying Lemma \ref{l51} we conclude that $\delta'$ is a univoque sequence in base $q$.
\smallskip

(b) Assume on the contrary that there exists a nontrivial univoque sequence  in some base $1<q\le p_m$. Then it is also univoque in base $p_m$. Furthermore, since a univoque sequence in a base $\le P_m$ contains at most finitely many zero digits and since a univoque sequence remains univoque if we remove an arbitrary finite initial block, the  there exists also a univoque sequence $(\eta_i)$ in base $p_m$ that contains only the digits $1$ and $m$.

It follows from the lexicographic characterization of univoque sequences that
\begin{equation*}
\eta_n=1\Longrightarrow (\lam_i)<(\eta_{n+i})<(\gamma_i)
\end{equation*}
and therefore (using the preceding lemma) that either
\begin{equation*}
\eta_n=1\Longrightarrow (\delta_{1+i})<(\eta_{n+i})<(\delta_i)
\end{equation*}
or
\begin{equation*}
\eta_n=1\Longrightarrow (\delta_i')<(\eta_{n+i})<m(\delta_i)
\end{equation*}

Setting $c_i=0$ if $\eta_i=1$ and $c_i=1$ if $\eta_i=m$ we obtain a sequence $(c_i)$ of zeroes and ones, satisfying either
\begin{equation}\label{511}
(d_{1+i})<(c_{n+i})<(d_i)\quad\text{whenever}\quad c_n=0
\end{equation}
or
\begin{equation}\label{512}
(d_i')<(c_{n+i})<1(d_i')\quad\text{whenever}\quad c_n=0.
\end{equation}
The second inequalities imply that $(c_i)$ has infinitely many zero digits. By removing a finite initial block if necessary we obtain a new sequence (still denoted by $(c_i)$) which begins with $c_1=0$ and which satisfies \eqref{511} or \eqref{512}.

In case of \eqref{511} we claim that
\begin{equation}\label{513}
0d_2d_3\cdots <(c_{n+i})<(d_i)\quad\text{for all}\quad n\ge 0.
\end{equation}
Indeed, if $c_n=1$ for some $n$ then there exist $m<n\le M$ such that $c_m=c_{M+1}=0$ and $c_{m+1}=\cdots =c_M=1$. In case of \eqref{511} it follows that
\begin{equation*}
(c_{n+i})\le (c_{m+i})<(d_i)
\end{equation*}
and
\begin{equation*}
(c_{n+i})\ge (c_{M+i})=0(c_{M+1+i})>0(d_{1+i})= 0d_2d_3\cdots .
\end{equation*}
However, \eqref{513} contradicts Lemma \ref{l35}.

In case of \eqref{512} we claim that
\begin{equation}\label{514}
0(d_i') <(c_{n+i})<1(d_i')\quad\text{for all}\quad n\ge 0.
\end{equation}
Indeed, if $c_n=1$ for some $n$ then there exist $m<n\le M$ such that $c_m=c_{M+1}=0$ and $c_{m+1}=\cdots =c_M=1$; then we have
\begin{equation*}
(c_{n+i})\le (c_{m+i})<1(d_i')
\end{equation*}
and
\begin{equation*}
(c_{n+i})\ge (c_{M+i})=0(c_{M+1+i})>0(d_i').
\end{equation*}
However, \eqref{514} contradicts Lemma \ref{l36}.
\end{proof}

The following lemma completes the proof of Theorem \ref{t11}.

\begin{lemma}\label{l59}\mbox{}

(a) If $d<\tilde d<1^{\infty}$ are admissible sequences, then $M_d\le m_{\tilde d}$ with equality if and only if $d=S(N,1)^{\infty}$ is of finite type and $\tilde d=S(N-1,1)S(N,1)^{\infty}$.
\smallskip

(b) The sets $I_d$, where $d$ runs over all admissible sequences $d\ne 1^{\infty}$, form a partition  of the interval $[ \frac{1+\sqrt{5}}{2},\infty )$.
\smallskip

(c) The set $C$ of numbers $m> \frac{1+\sqrt{5}}{2}$ satisfying $p_m=P_m$ is a Cantor set, i.e., a nonempty closed set having neither interior, nor isolated points. Its smallest element is $1+x\approx 2.3247$ where $x$ is the first Pisot number, i.e., the positive root of the equation $x^3=x+1$.
\end{lemma}

\begin{proof}\mbox{}

(a) It suffices to prove that if $d<\tilde d<1^{\infty}$ are admissible sequences \emph{of infinite type}, then $M_d< M_{\tilde d}$. (We recall that $m_d=M_d$ and $m_{\tilde d}=M_{\tilde d}$ in this case.) Indeed, the general case hence follows by recalling from Example \ref{e32} and Lemma \ref{l56} that if $d=S(N,1)^{\infty}$ is of finite type, then $\tilde d=S(N-1,1)S(N,1)^{\infty}$ is the smallest admissible sequence satisfying $\tilde d>d$, and $m_d<M_d=m_{\tilde d}= M_{\tilde d}$.

Furthermore, it is sufficient to show that if $d<\tilde d<1^{\infty}$ are admissible sequences \emph{of infinite type}, then $p_{d,m}''>p_{\tilde d,m}''$ for each $m\in (1,\infty)$ where $p_{d,m}''$ and $p_{\tilde d,m}''$ denote the expressions $p_m''$ of Section \ref{s4} for the admissible sequences $d$ and $\tilde d$, respectively. Indeed, then  we can conclude that $p_{d,M_{\tilde d}}''>p_{\tilde d,M_{\tilde d}}''=P_{M_{\tilde d}}$ and therefore, since the function $m\mapsto p_{d,m}''-P_m$ is strictly increasing Lemma \ref{l41}, $M_d< M_{\tilde d}$.

Assuming on the contrary that $p_{d,m}''\le p_{\tilde d,m}''$ for some $m$, in base $q:=p_{\tilde d,m}''$ we have
\begin{equation*}
\pi_q(m-\tilde \delta')=1=\pi_{p_{d,m}''}(m- \delta')\ge \pi_q(m- \delta')
\Longrightarrow \pi_q(\delta')\ge \pi_q(\tilde \delta')
\end{equation*}
Since $d$ and $\tilde d$ are of infinite type, we have $\delta=m\delta'$ and $\tilde \delta=m\tilde \delta'$ by Lemma \ref{l34}, so that the last inequality is equivalent to $\pi_q(\delta)\ge \pi_q(\tilde \delta)$.

Since  quasi-greedy expansions remain quasi-greedy in larger bases, it follows from Lemma \ref{l56} that both $\delta$ and $\tilde \delta$ are quasi-greedy expansions in base $q$. Therefore we deduce from the last inequality that $\delta\ge\tilde \delta$, contradicting our assumption.
\smallskip

(b) The sets $I_d$ are disjoint by (a) and they cover the interval $\left[ \frac{1+\sqrt{5}}{2},\infty \right) $ by Lemma \ref{l54}. In view of (a) the proof will be completed if we show that for the smallest admissible sequence we have
\begin{equation}\label{515}
I_{0^{\infty}}=\left[ \frac{1+\sqrt{5}}{2},1+P_1\right)
\end{equation}
where $x>1$ is the first Pisot number.

The values $m_d$ and $M_d$ are the solutions of the equations
\begin{equation*}
\pi_{P_m}(\delta)=m-1
\quad\text{and}\quad
\pi_{P_m}(\delta')=\frac{m}{P_m-1}-1.
\end{equation*}
Now we have $\delta=\delta'=1^{\infty}$, so that our equations take the form
\begin{align*}
&\frac{1}{P_m-1}=m-1
\intertext{and}
&\frac{1}{P_m-1}=\frac{m}{P_m-1}-1.
\end{align*}
Using \eqref{43} we obtain that they are equivalent to $m=(1+\sqrt{5})/2$ and $m=1+P_1$, respectively.
\smallskip

(c) If we denote by $D_1$ and $D_2$ the set of admissible sequences $d\ne 1^{\infty}$ of finite and infinite type, respectively, then
\begin{equation*}
C=[2,\infty)\setminus\cup_{d\in D_1}(m_d,M_d)
\end{equation*}
so that $C$ is a closed set. The relation \eqref{515} shows that its smallest element is $1+P_1$. In order to prove that it is a Cantor set, it suffices to show that
\begin{itemize}
\item the intervals $[m_d,M_d]$ ($d\in D_1)$ are disjoint;

\item for each $m\in C$ there exist two sequences $(a_N)\subset [2,\infty)\setminus C$ and $(b_N)\subset C\setminus \set{m}$, both converging to $m$.
\end{itemize}

The first property follows from (a). For the proof of the second property let us consider the infinite sequence $h=(h_j)$ of positive integers defining the admissible sequence $d$ for which $m_d=m$, and  set $d_N:=S_h(N,1)^{\infty}$, $N=1,2,\ldots .$ This is a decreasing sequence of admissible sequences, converging pointwise to $d$. Using (a) we conclude that both $(m_{d_N})$ and $(M_{d_N})$ converge to $m_d=M_d$. Since $m_{d_N}\in D_1$ and $M_{d_N}\in D_2$ for every $N$, the proof is complete.
\end{proof}

\end{document}